\documentclass[12pt]{article}
\usepackage{amsmath,amssymb,amsbsy,amsfonts,amsthm,
     latexsym,amsopn,amstext,amsxtra,euscript,amscd}

\begin{document}

\newtheorem{lem}{Lemma}
\newtheorem{lemma}[lem]{Lemma}
\newtheorem{prop}{Proposition}
\newtheorem{proposition}[prop]{Proposition}
\newtheorem{thm}{Theorem}
\newtheorem{theorem}[thm]{Theorem}
\newtheorem{cor}{Corollary}
\newtheorem{corollary}[cor]{Corollary}
\newtheorem{prob}{Problem}
\newtheorem{problem}[prob]{Problem}
\newtheorem{ques}{Question}
\newtheorem{question}[ques]{Question}
\newtheorem{rems}{Remarks}
\newtheorem{remarks}[rems]{Remarks}

\def\mand{\qquad\mbox{and}\qquad}
\def\scr{\scriptstyle}
\def\\{\cr}
\def\({\left(}
\def\){\right)}
\def\[{\left[}
\def\]{\right]}
\def\<{\langle}
\def\>{\rangle}
\def\fl#1{\left\lfloor#1\right\rfloor}
\def\rf#1{\left\lceil#1\right\rceil}
\def\cA{{\mathcal A}}
\def\cB{{\mathcal B}}
\def\cC{{\mathcal C}}
\def\cE{{\mathcal E}}
\def\cI{{\mathcal I}}
\def\cJ{{\mathcal J}}
\def\cK{{\mathcal K}}
\def\cN{{\mathcal N}}
\def\cM{{\mathcal M}}
\def\cP{{\mathcal P}}
\def\cQ{{\mathcal Q}}
\def\cR{{\mathcal R}}
\def\cS{{\mathcal S}}
\def\cT{{\mathcal T}}
\def\N{{\mathbb N}}
\def\Z{{\mathbb Z}}
\def\eps{\varepsilon}

\newcommand{\comm}[1]{\marginpar{\fbox{#1}}}
\def\xxx{\vskip5pt\hrule\vskip5pt}
\def\yyy{\vskip5pt\hrule\vskip2pt\hrule\vskip5pt}

\title{\bf Noncototients and Nonaliquots}

\author{
{\sc William D.~Banks} \\
{Department of Mathematics, University of Missouri} \\
{Columbia, MO 65211, USA} \\
{\tt bbanks@math.missouri.edu} \\
\and
{\sc Florian~Luca} \\
{Instituto de Matem{\'a}ticas}\\
{ Universidad Nacional Aut\'onoma de M{\'e}xico} \\
{C.P. 58180, Morelia, Michoac{\'a}n, M{\'e}xico} \\
{\tt fluca@matmor.unam.mx}}

\date{\today}

\maketitle

\begin{abstract}
Let $\varphi(\cdot)$ and $\sigma(\cdot)$ denote the Euler
function and the sum of divisors function, respectively.
In this paper, we give a lower bound for the number of positive
integers $m\le x$ for which the equation $m=n-\varphi(n)$ has no solution.
We also give a lower bound for the number of $m\le x$ for which
the equation $m=\sigma(n)-n$ has no solution.  Finally, we show the set
of positive integers $m$ not of the form $(p-1)/2-\varphi(p-1)$
for some prime number $p$ has a positive lower asymptotic density.
\end{abstract}

\newpage

\section{Introduction}

Let $\varphi(\cdot)$ denote the Euler function, whose value at the
positive integer $n$ is
$$
\varphi(n)=n\prod_{p|n}\(1-\frac{1}{p}\).
$$
An integer of the form $\varphi(n)$ is called a
{\it totient\/}; a {\it cototient\/} is an integer in the
image of the function $f_c(n)=n-\varphi(n)$.  If $m$ is a positive
integer for which the equation $f_c(n)=m$ has no solution, then $m$
is called a {\it noncototient\/}.  An old conjecture of
Erd\H os and Sierpi\'nski (see B36 in~\cite{Guy}) asserts the
existence of infinitely many noncototients.
This conjecture has been settled by Browkin and Schinzel~\cite{BS},
who showed that if $w\ge 3$ is an odd integer satisfying certain
arithmetic properties, then $m=2^\ell w$ is a noncototient for every
positive integer $\ell$; they also showed that the integer $w=509203$
is one such integer.  Flammenkamp and Luca~\cite{FL} later found
six more integers $w$ satisfying the same properties.
These results, however, imply only the weak lower bound
$\#\cN_c(x)\gg\log x$ for the cardinality of the set
$$
\cN_c(x)=\{1\le m\le x:
m\ne f_c(n)\text{~for every positive integer~}n\}.
$$
In Theorem~\ref{thm:main1} (Section~\ref{sec:noncototient}),
we show that $2p$ is a noncototient for almost every
prime $p$ (that is, for all $p$ in a set of primes of relative
asymptotic density one), which implies the following unconditional
lower bound for the number of noncototients $m\le x$:
$$
\#\cN_c(x)\ge\frac{x}{2\log x}\,(1+o(1)).
$$

Next, let $\sigma(\cdot)$ denote the sum of divisors function,
whose value at the positive integer $n$ is
$$
\sigma(n)=\sum_{d|n}d=\prod_{p^a\|n}\frac{p^{a+1}-1}{p-1}.
$$
An integer in the image of the function $f_a(n)=\sigma(n)-n$
is called an {\it aliquot number\/}. If $m$ is a positive
integer for which the equation $f_a(n)=m$ has no solution, then
$m$ is said to be {\it nonaliquot\/}.
Erd\H os~\cite{Er} showed that the collection of nonaliquot
numbers has a positive lower asymptotic density,
but no numerical lower bound on this density was given.
In Theorem~\ref{thm:main2} (Section~\ref{sec:nonaliquot}),
we show that the lower bound
$\#\cN_a(x)\ge\tfrac{1}{48}x\,(1+o(1))$ holds, where
$$
\cN_a(x)=\{1\le m\le x:
m\ne f_a(n)\text{~for every positive integer~}n\}.
$$

Finally, for an odd prime $p$, let $f_r(p)=(p-1)/2-\varphi(p-1)$.
Note that $f_r(p)$ counts the number of quadratic nonresidues
modulo~$p$ which are not primitive roots. 
At the 2002 Western Number Theory Conference in San Francisco,
Neville Robbins asked whether there exist infinitely many 
positive integers $m$ for which $f_r(p)=m$ has no solution;
let us refer to such integers as {\it Robbins numbers\/}.
The existence of infinitely many Robbins numbers has
been shown recently by Luca and Walsh~\cite{LuWa}, who proved
that for every odd integer $w\ge 3$, there exist infinitely 
many integers $\ell\ge 1$ such that $2^\ell w$ is a Robbins number.
In Theorem~\ref{thm:main3} (Section~\ref{sec:Robbins}),
we show that the set of Robbins numbers has a positive density;
more precisely, if
$$
\cN_r(x)=\{1\le m\le x:m\ne f_r(p)\text{~for every odd prime~}p\},
$$ 
then the lower bound $\#\cN_r(x)\ge\tfrac{1}{3}x\,(1+o(1))$ holds.

\bigskip

{\bf Notation.}
Throughout the paper, the letters $p$, $q$ and $r$ are always used
to denote prime numbers.  For an integer $n\ge 2$, we write $P(n)$
for the largest prime factor of $n$, and we put $P(1)=1$.
As usual, $\pi(x)$ denotes the number of primes 
$p\le x$, and if $a,b>0$ are coprime integers, 
$\pi(x;b,a)$ denotes the number of primes $p\le x$ such that
$p\equiv a\pmod b$.  For any set $\cA$ and real number $x\ge 1$,
we denote by $\cA(x)$ the set $\cA\cap[1,x]$.
For a positive integer $k$, we write $\log_k(\cdot)$
for the function given recursively by $\log_1 x=\max\{\log x,1\}$
and $\log_k x=\log_1(\log_{k-1}x)$, where $x>0$ is a real number
and $\log(\cdot)$ denotes the natural logarithm. When $k=1$,
we omit the subscript in order to simplify the notation, with the
continued understanding that $\log x\ge 1$ for all $x>0$.
We use the Vinogradov symbols $\ll$ and $\gg$,
as well as the Landau symbols $O$ and $o$, with their usual meanings. 
Finally, we use $c_1,c_2,\ldots$ to denote constants that
are positive and absolute.

\bigskip

{\bf Acknowledgements.}
Most of this work was done during a visit by the second author
to the University of Missouri--Columbia; the hospitality and
support of this institution are gratefully acknowledged.
During the preparation of this paper,
W.~B.\ was supported in part by NSF grant DMS-0070628, and 
F.~L.\ was supported in part by grants
SEP-CONACYT 37259-E and 37260-E.

\section{Noncototients}
\label{sec:noncototient}

We begin this section with some technical results that are needed for the
proof of Theorem~\ref{thm:main1} below.

\begin{lem}
\label{lem:1}
The following estimate holds:
$$
\sum_{\substack{x^{1-1/t}<n\le x^{1-1/(t+1)}\\P(n)\le x/n}}
\frac{1}{n}\ll \log x\left\{\begin{array}{ll}
\exp(-0.5t\log t)&\quad\text{if $t\le(\log x)/(3\log_2x)$};\\
\exp(-0.5t)&\quad\text{otherwise}.\end{array}\right.
$$
\end{lem}

\begin{proof}
For all $x\ge y\ge 2$, let
$$
\Psi(x,y)=\#\{n\le x:P(n)\le y\},
$$
and put $u=(\log x)/(\log y)$.  If $u\le y^{1/2}$, the estimate
\begin{equation}
\label{eq:psi1}
\Psi(x,y)=xu^{-u+o(u)}
\end{equation}
holds (see Corollary~1.3 of~\cite{HT}, or~\cite{CEP}), while the
upper bound
\begin{equation}
\label{eq:psi2}
\Psi(x,y)\ll xe^{-u/2}
\end{equation}
holds for arbitrary $u\ge 1$ (see, for example, Theorem~1 in
Chapter~III.5 of~\cite{Tenen}).  Since
$$
\sum_{\substack{x^{1-1/t}<n\le x^{1-1/(t+1)}\\P(n)\le x/n}}1
\le\Psi\big(x^{1-1/(t+1)},x^{1/t}\big),
$$
the result follows from~\eqref{eq:psi1} and~\eqref{eq:psi2} by
partial summation.
\end{proof}

For every integer $n\ge 3$ and real number $y>2$, let
$$
h_y(n)=\sum_{\substack{p|(2n-\varphi(n))\\ p>y}}\frac{1}{p}.
$$

\begin{lem}
\label{lem:2}
Let $\cA$ be the set of integers $n\ge 3$ for which
$\gcd(n,\varphi(n))=1$, and let
$$
\cA(x,y)=\{n\in\cA(x):h_y(n)>1\}.
$$
Then, uniformly for $2<y\le(\log x)^{1/4}$, the following
estimate holds: 
$$
\sum_{n\in \cA(x,y)}\frac{1}{n}\ll\frac{\log x}{y\log_2y}.
$$
\end{lem}

\begin{proof} Our proof follows closely 
the proof of Lemma 3 from \cite{LuPo}.

We first determine an upper bound 
on the cardinality $\#\cA(x,y)$ of the set $\cA(x,y)$
in the case that $2<y\le(\log x)^{1/2}$.  Let
$$
z=\exp\(\frac{\log x\log_2 y}{2\log y}\)\qquad\text{and}\qquad
u=\frac{\log x}{\log z}=\frac{2\log y}{\log_2y}.
$$
Then
$$
u\log u=2(1+o(1))\log y.
$$
Let $\cA_1(x,y)=\{n\in\cA(x):P(n)\le z\}$.  Since $y\le(\log x)^{1/2}$,
it follows that $u\le z^{1/2}$; therefore, using~\eqref{eq:psi1} we
derive that
\begin{equation}
\label{eq:A1}
\#\cA_1(x,y)\le \Psi(x,z)=\frac{x}{\exp((1+o(1))u\log u)}=
\frac{x}{y^{2+o(1)}}\ll \frac{x}{y\log_2y}.
\end{equation}
For each $n\in \cA(x,y)\backslash \cA_1(x,y)$, write $n$ in the form 
$n=Pk$, where $P>z$ is prime, and $k<x/z$.  Note that $n$ is
squarefree since $\gcd(n,\varphi(n))=1$.
Let $\cA_2(x,y)$ be the set of those integers
$n\in\cA(x,y)\backslash\cA_1(x,y)$ for which $k\le 2$.
Clearly,
\begin{equation}
\label{eq:A2}
\#\cA_2(x,y)\le \pi(x)+\pi(x/2)
\ll\frac{x}{\log x}\le\frac{x}{y\log_2y}.
\end{equation}
Now let $\cA_3(x,y)=\cA(x,y)\backslash\(\cA_1(x,y)\cup\cA_2(x,y)\)$,
and suppose that $n$ lies in $\cA_3(x,y)$.
For a fixed prime $p>y$, if $p|(2n-\varphi(n))$, then
\begin{equation}
\label{eq:P}
P(2k-\varphi(k))+\varphi(k)\equiv 0\pmod p.
\end{equation}
Fixing $k$ as well, we see that $p\ne P$ 
(otherwise, $P|\varphi(k)|\varphi(n)$ and $P|n$, which contradicts
the fact that $n\in\cA$), 
and $p\nmid(2k-\varphi(k))$ (otherwise, it follows that
$p|\gcd(k,\varphi(k))|\gcd(n,\varphi(n))=1$).
Let $a_k$ be the congruence class 
modulo $p$ determined for $P$ by the congruence~\eqref{eq:P}; 
then the number of possibities for~$n$ (with $p$ and $k$ fixed)
is at most $\pi(x/k;p,a_k)$.

In the case that $pk\le x/z^{1/2}$, we use a well known result
of Montgomery and Vaughan~\cite{MV} to conclude that 
$$
\pi(x/k;p,a_k)\le\frac{2x}{\varphi(p)k\log(x/kp)}\le 
\frac{4x}{(p-1)k\log z}
\le\frac{12x\log y}{pk\log x\log_2y}.
$$
In the case that $x/z^{1/2}<pk<x$, since $k<x/z$, we
see that $p>z^{1/2}$.  Here, we use the trivial estimate 
$$
\pi(x/k;p,a_k)\le\frac{x}{pk}.
$$
Finally, if $pk\ge x$, then $p>z$, and we have 
$$
\pi(x/k;p,a_k)\le 1.
$$

Now, for fixed $p>y$, let
$$
\cA_3(p,x,y)=\{n\in\cA_3(x,y):p|(2n-\varphi(n))\}.
$$ 
When $p\le z^{1/2}$, we have
$$
\#\cA_3(p,x,y)\le 
\frac{12x\log y}{p\log x\log_2y}\sum_{k<x/z}\frac{1}{k}\ll 
\frac{x\log y}{p\log_2y}.
$$
If $z^{1/2}<p\le z$, then
\begin{eqnarray*}
\#\cA_3(p,x,y)&\le&
\frac{12x\log y}{p\log x\log_2y}
\sum_{k<x/z}\frac{1}{k}+\frac{x}{p}\sum_{k<x/z}\frac{1}{k}\\
&\ll&\frac{x\log y}{p\log_2y}+\frac{x\log x}{p}\ll\frac{x\log x}{p}.
\end{eqnarray*}
Finally, if $p>z$, it follows that 
\begin{eqnarray*}
\#\cA_3(p,x,y)&\le&\frac{12x\log y}{p\log x\log_2y}
\sum_{k<x/z}\frac{1}{k}+\frac{x}{p}\sum_{k<x/z}\frac{1}{k}+
\sum_{k<x/z} 1\\
&\ll&\frac{x\log y}{p\log_2y}+\frac{x\log x}{p}+\frac{x}{z}
\ll\frac{x\log x}{z}.
\end{eqnarray*}
Consequently,
\begin{eqnarray}
\#\cA_3(x,y)&=&\sum_{n\in\cA_3(x,y)} 1
<\sum_{n\in\cA_3(x,y)}h_y(n)\nonumber\\
&=&\sum_{n\in\cA_3(x,y)}
\sum_{\substack{p|(2n-\varphi(n))\\ p>y}}\frac{1}{p}
=\sum_{p>y}\frac{1}{p}\#\cA_3(p,x,y)\nonumber\\
&\ll&\frac{x\log y}{\log_2y}
\sum_{y<p\le z^{1/2}}\frac{1}{p^2}
+x\log x\sum_{z^{1/2}<p\le z}\frac{1}{p^2}
+\frac{x\log x}{z}\sum_{z<p\le 2x}\frac{1}{p}\nonumber\\
\label{eq:A3}
&\ll&\frac{x}{y\log_2y}
+\frac{x\log x}{z^{1/2}\log z}
+\frac{x\log x\log_2x}{z}\ll\frac{x}{y\log_2y},
\end{eqnarray}
where the last estimates follows (if $x$ is sufficiently large)
from the bound $y\le(\log x)^{1/2}$ and our choice of $z$.
Thus, by the inequalities~\eqref{eq:A1}, \eqref{eq:A2},
and~\eqref{eq:A3}, we obtain that
$$
\#\cA(x,y)\ll\frac{x}{y\log_2y}.
$$

Now, for all $y\le(\log x)^{1/4}$, we have by partial summation
(using the fact that $y\le(\log t)^{1/2}$ if $\exp(y^2)\le t\le x$):
\begin{eqnarray*}
\sum_{n\in\cA(x,y)}\frac{1}{n}
&\le&\sum_{n\le\exp(y^2)}\frac{1}{n}
+\sum_{\substack{\exp(y^2)\le n\le x\\
n\in \cA(x,y)}}\frac{1}{n}\\
&\ll& y^2+\int_{\exp(y^2)}^x\frac{d\cA(t,y)}{t}\\
&\ll& y^2+\frac{1}{y\log_2y}\int_{1}^x\frac{dt}{t}
=y^2+\frac{\log x}{y\log_2y}
\ll\frac{\log x}{y\log_2y},
\end{eqnarray*}
which completes the proof.
\end{proof}

\begin{lem}
\label{lem:3}
For some absolute constant $c_1>0$, the set $\cB$ defined by
$$
\cB=\big\{n:p\nmid\varphi(n)
\text{~for some prime~}p\le c_1(\log_2n)/(\log_3n)\big\}
$$
satisfies
$$
\sum_{n\in\cB(x)}\frac{1}{n}\ll\frac{\log x}{(\log_2 x)^2}.
$$
\end{lem}

\begin{proof}
By Theorem 3.4 in \cite{EGPS}, there exist positive constants 
$c_0,c_2,x_0$ such that for all $x\ge x_0$, the bound
$$
S'(x,p)=\sum_{\substack{q\le x\\p|(q-1)}}\!\!\!\vphantom{\biggl|}'
~\frac{1}{q}\ge\frac{c_2\log_2 x}{p},
$$
where the dash indicates that the prime $q$ is omitted from the sum
if there exists a real primitive character $\chi$ modulo $q$
for which $L(s,\chi)$ has a real root $\beta\ge 1-c_0/\log q$.
{From} the proof of Theorem~4.1 in~\cite{EGPS}, we also have
the estimate
$$
\sum_{\substack{n\le x\\p\,\nmid\,\varphi(n)}} 1
\ll \frac{x}{\exp(S'(x,p))},
$$
uniformly in $p$. Thus, if $c_1=c_2/3$, $g(x)=c_1(\log_2 x)/(\log_3 x)$,
and $p\le g(x)$, then 
$$
\sum_{\substack{n\le x\\ p\,\nmid\,\varphi(n)}} 1
\ll \frac{x}{(\log_2 x)^3}.
$$
Therefore,
$$
\sum_{p\le g(x)}\sum_{\substack{n\le x\\ p\,\nmid\,\varphi(n)}} 1\ll 
\frac{x\pi(g(x))}{(\log_2 x)^3}\ll \frac{x}{(\log_2 x)^2}.
$$
This argument shows that the inequality 
$$
\#\cB(x)\ll \frac{x}{(\log_2 x)^2}
$$
holds uniformly in $x$, and the result follows by partial summation. 
\end{proof}

The following lemma is a consequence of well known
estimates for the number of integers $n\le x$ free of
prime factors $p\le y$. In particular, the result follows
immediately, using partial summation, from Theorem~3 and
Corollary~3.1 in Chapter~III.6 of~\cite{Tenen}; the proof
is omitted.

\begin{lem}
\label{lem:4}
Let
$$
\cC(x;y)=\{n\le x:p\nmid n\text{~for all~}p\le y\}.
$$
Then, uniformly for $2\le y\le(\log x)^{1/2}$, we have 
$$
\sum_{n\in\cC(x;y)}\frac{1}{n}\ll\frac{\log x}{\log y}.
$$
\end{lem}

We now come to the main result of this section.

\begin{thm}
\label{thm:main1}
For almost all primes $p$ (that is, for all primes $p$ in a set
of relative asymptotic density $1$), the number $2p$ is a
noncototient: $2p\in\cN_c$.
In particular, the inequality 
$$
\#\cN_c(x)\ge\frac{x}{2\log x}(1+o(1))
$$
holds as $x\rightarrow\infty$.
\end{thm}

\begin{proof}
Suppose that
$$
f_c(n)=n-\varphi(n)=2p
$$
holds, where $p\le x/2$ is an odd prime. We can assume that
$p>x/\log x$, since the number of primes $p\le x/\log x$ is 
$\pi(x/\log x)=o(\pi(x/2))$. Then $n\ge 3$, and $\varphi(n)$ is
even; hence, $n$ is also even. If $4|n$, then $2\|\varphi(n)$,
and the only possibility is $n=4$, which is not possible.
Thus, $2\|n$.  Writing $n=2m$, with $m$ odd, the
equation above becomes
\begin{equation}
\label{eq:2m}
f_c(2m)=2m-\varphi(m)=2p.
\end{equation}
Clearly, $x\ge 2p\ge 2m-\varphi(m)\ge m$. Now observe
that $\gcd(m,\varphi(m))=1$. Indeed, if $q|\gcd(m,\varphi(m))$
for an odd prime $q$, it must be the case that $q=p$. 
Then, either $p^2|m$, or $pr|m$ for some prime $r\equiv 1\pmod p$. 
In both cases, we see that $x\ge m\ge p^2\ge (x/\log x)^2$,
which is not possible since $m\le x$.
In particular, $m$ lies in the set $\cA(x)$ defined in Lemma~\ref{lem:2}.
Finally, we can assume that $m$ is not prime, for
otherwise~\eqref{eq:2m} becomes $m=2p-1$, which is well known
to have at most $O(x/(\log x)^2)=o(\pi(x/2))$ solutions with
primes $m,p$ such that $p\le x/2$.

Let $\cM(x)$ be the set of (squarefree odd)
integers $m$ for which~\eqref{eq:2m} holds for some
prime $p>x/\log x$. To prove the theorem, it suffices to
show that $\#\cM(x)=o(x/\log x)$.

Let $m\in\cM(x)$, and write $m=Pk$, where $P=P(m)>P(k)$
and $k\ge 3$. Since $m>p>x/\log x$ is squarefree, it 
follows that $P\gg\log x$. Equation~\eqref{eq:2m} now becomes 
$$
P(k-\varphi(k)/2)-\varphi(k)=p.
$$
For fixed $k$, we apply the sieve (see, for example,
Theorem~5.7 of~\cite{HR}) to conclude that 
the number of possibilities for $P$ (or $p$) is 
\begin{eqnarray}
\lefteqn{\ll \frac{x}{\varphi(k-\varphi(k)/2)}\cdot 
\frac{1}{\(\log\(x/(k-\varphi(k)/2)\)\)^2}}\nonumber\\
\label{eq:choiceforP}
&&\qquad\qquad\qquad\ll\frac{x}{\varphi(k-\varphi(k)/2)}
\cdot\frac{1}{\(\log(x/k)\)^2}.
\end{eqnarray}

Now put 
$$
y_1=\exp\(\frac{\log x\log_4 x}{3\log_3 x}\),
$$
and let $\cM_1(x)=\{m\in\cM(x):P\le y_1\}$. 
For $m\in\cM_1(x)$, we have
$$
k>\frac{x}{P\log x}\ge\frac{x}{y_1\log x}.
$$
In particular, if $x$ is sufficiently large, and
$t_1=4(\log_3 x)/(\log_4 x)$, then every integer $k$ belongs
to an interval of the form $\cI_j=[x^{1-1/(t_1+j)},x^{1-1/(t_1+j+1)}]$
for some nonnegative integer $j$ such that $t_1+j+1\le \log x$.
For fixed $j$, we have $\log(x/k)\gg (\log x)/(t_1+j)$, and therefore 
$$
\frac{1}{(\log(x/k))^2}\ll \frac{(t_1+j)^2}{(\log x)^2}.
$$
Using the fact that $\varphi(n)\gg n/\log_2 n$, we see that
for each fixed $k\in\cI_j$, the number of choices for $P$ is
$$
\ll \frac{x\log_2 x}{(\log x)^2}\cdot \frac{(t_1+j)^2}{2k-\varphi(k)}
<\frac{x\log_2 x}{(\log x)^2}\cdot \frac{(t_1+j)^2}{k}.
$$
Summing first over $k$, then $j$, and applying Lemma~\ref{lem:1}, 
we derive that
\begin{eqnarray}
\label{eq:M1}
\#\cM_1(x)&\ll&\frac{x\log_2 x}{(\log x)^2}
\sum_{0\le j\le \log x-t_1}(t_1+j)^2
\sum_{\substack{k\in \cI_j\\P(k)<x/k}}\frac{1}{k}
\nonumber\\
&\ll&\frac{x\log_2 x}{\log x}
\sum_{0\le j\le(\log x)/(3\log_2x)-t_1}
\frac{(t_1+j)^2}{\exp\(0.5(t_1+j)\log(t_1+j)\)}\nonumber\\
&&\qquad+\quad\frac{x\log_2 x}{\log x}
\sum_{j>(\log x)/(3\log_2x)-t_1}
\frac{(t_1+j)^2}{\exp\(0.5(t_1+j)\)}
\nonumber\\
&\ll&\frac{x\log_2 x}{\log x}\cdot
\frac{t_1^2}{\exp(0.5t_1\log t_1)}
+\frac{x\log x}{\log_2x}\exp\(-\frac{\log x}{6\log_2x}\)
\nonumber\\
&\ll&
\frac{x\log_2 x}{\log x}\frac{(\log_3x)^2}{(\log_4x)^2}
\exp\(-2(1+o(1))\log_3x\)
+o\(\frac{x}{\log x}\)\nonumber \\
& = & o\(\frac{x}{\log x}\).
\end{eqnarray}

Hence, from now on, we need only consider numbers
$m\in\cM(x)\backslash\cM_1(x)$.  For such integers, we have 
$x/k\ge P>y_1$; thus,
$$
\frac{1}{\(\log(x/k)\)^2}\ll \frac{1}{(\log y_1)^2}\ll 
\frac{(\log_3 x)^2}{(\log x\log_4 x)^2}.
$$
For fixed $k$, the number of choices~\eqref{eq:choiceforP} 
for the prime $P$ is 
$$
\ll \frac{x(\log_3 x)^2}{(\log x\log_4 x)^2}\cdot \frac{1}
{\varphi(k-\varphi(k/2))}.
$$
Put $y_2=\exp\(\exp\(\sqrt{\log_3 x}\,\)\)$, and let 
$$
\cM_2(x)=\{m\in \cM(x)\backslash \cM_1(x):k\in\cA(x,y_2)\},
$$
where $\cA(x,y_2)$ is defined as in Lemma~\ref{lem:2}.
Using once more the inequality 
$\varphi(n)\gg n/\log_2 n$, the fact that $k-\varphi(k)/2\ge k/2$,
and Lemma~\ref{lem:2}, we have
\begin{eqnarray}
\label{eq:M2}
\#\cM_2(x)& \ll & \frac{x(\log_3 x)^2\log_2 x}{(\log x\log_4 x)^2}
\sum_{k\in\cA(x,y_2)}\frac{1}{k}\nonumber \\
&\ll&\frac{x(\log_2 x)^2}{y_2\log x}
=o\(\frac{x}{\log x}\)
\end{eqnarray}
since $(\log_2 x)^2=o(y_2)$.

Next, we consider numbers $m\in\cM(x)$ that do not
lie in $\cM_1(x)\cup\cM_2(x)$. For such integers, we have
$$
\sum_{p|(2k-\varphi(k))}\frac{1}{p}
\le\sum_{p\le y_2}\frac{1}{p}+1=
\log_2 y_2+O(1)={\sqrt {\log_3 x}}+O(1).
$$
Therefore,
\begin{eqnarray*}
\frac{1}{\varphi(k-\varphi(k)/2)}&=&\frac{1}{k-\varphi(k)/2}\cdot 
\frac{k-\varphi(k/2)}{\varphi(k-\varphi(k)/2)}\\
& \ll & 
\frac{1}{k}\prod_{p|(2k-\varphi(k))}\(1+\frac{1}{p-1}\)\le 
\frac{1}{k}\exp\(\sum_{p|(2k-\varphi(k))}\frac{1}{p}\)\\
& \ll & 
\frac{\exp\(\sqrt {\log_3 x}\,\)}{k}.
\end{eqnarray*}
Now put 
$$
y_3=\exp\(\log x\(\frac{\log_4x}{\log_3x}\)^{1/2}\),
$$
and let 
$$
\cM_3(x)=\big\{m\in \cM(x)\backslash\(\cM_1(x)\cup\cM_2(x)\):
P(m)\le y_3\big\}.
$$ 
In this case, 
$$
k>\frac{x}{P\log x}>\frac{x}{y_3\log x}.
$$
In particular, if $x$ is sufficiently large, and
$t_2=2((\log_3 x)/(\log_4 x))^{1/2}$, every such $k$ belongs
to an interval of the form $\cJ_j=[x^{1-1/(t_2+j)},x^{1-1/(t_2+j+1)}]$
for some nonnegative integer $j$ such that $t_2+j+1\le \log x$.
For fixed $j$, we have
$\log(x/k)\gg (\log x)/(t_2+j)$, and therefore 
$$
\frac{1}{(\log(x/k))^2}\ll \frac{(t_2+j)^2}{(\log x)^2}.
$$
Using the fact that $\varphi(n)\gg n/{\exp\(\sqrt{\log_3 x}\,\)}$
for $n=k-\varphi(k)/2$, it follows that 
for any fixed $k\in \cJ_j$, the number of choices for $P$ is
$$
\ll \frac{x\exp\(\sqrt{\log_3 x}\,\)}
{(\log x)^2}\cdot \frac{(t_2+j)^2}{2k-\varphi(k)}
\ll \frac{x\exp\(\sqrt{\log_3 x}\,\)}
{(\log x)^2}\cdot \frac{(t_2+j)^2}{k}.
$$
Summing up first over $k$, then over $j$, and using Lemma~\ref{lem:1}
again, we obtain that 
\begin{eqnarray}
\label{eq:M3}
\lefteqn{\#\cM_3(x)
\ll\frac{x\exp\(\sqrt{\log_3 x}\,\)}
{(\log x)^2}\sum_{0\le j\le \log x-t_2}(t_2+j)^2
\sum_{\substack{k\in \cI_j\\P(k)<x/k}}\frac{1}{k}}\nonumber\\
&&\qquad\ll\frac{x\exp\(\sqrt{\log_3 x}\,\)}{\log x}
\sum_{0\le j\le(\log x)/(3\log_2x)-t_2}
\frac{(t_2+j)^2}{\exp\(0.5(t_2+j)\log(t_2+j)\)}\nonumber\\
&&\qquad\qquad+\quad\frac{x\exp\(\sqrt{\log_3 x}\,\)}{\log x}
\sum_{j>(\log x)/(3\log_2x)-t_2}
\frac{(t_2+j)^2}{\exp\(0.5(t_2+j)\)}\nonumber\\
&&\qquad\ll\frac{x\exp\(\sqrt{\log_3 x}\,\)}{\log x}\cdot
\(\frac{t_2^2}{\exp(0.5t_2\log t_2)}+
\exp\(-\frac{\log x}{6\log_2x}\)\)\nonumber\\
&&\qquad=o\(\frac{x}{\log x}\).
\end{eqnarray}

Hence, we can now restrict our attention to numbers $m\in\cM(x)$ 
which do not lie in $\cup_{i=1}^3\cM_i(x)$. 
For such numbers, we have $x/k\ge P>y_3$; thus, 
$$
\frac{1}{\(\log(x/k)\)^2}\ll \frac{1}{(\log y_3)^2}\ll 
\frac{\log_3 x}{(\log x)^2\log_4 x},
$$
and the number of choices~\eqref{eq:choiceforP} for $P$,
for fixed $k$, is 
\begin{equation}
\label{eq:newchoiceforP}
\ll \frac{x\log_3 x}{(\log x)^2\log_4 x}\cdot \frac{1}
{\varphi(k-\varphi(k/2))}.
\end{equation}
Let 
$$
\cM_4(x)=\left\{m\in\cM(x)\backslash\(\cup_{i=1}^3\cM_i(x)\):
k\le\exp\(\sqrt{\log x}\,\)\right\}.
$$ 
Clearly, by~\eqref{eq:newchoiceforP}, we have
\begin{eqnarray}
\label{eq:M4}
\#\cM_4(x)& \ll & \frac{x\log_2 x\log_3 x}{(\log x)^2\log_4 x}
\sum_{k\le\exp(\sqrt{\log x}\,)}\frac {1}{k}\nonumber \\
& \ll & 
\frac{x\log_2 x\log_3 x}{(\log x)^{3/2}\log_4 x}
=o\(\frac{x}{\log x}\).
\end{eqnarray}

Now let $\cB$ be the set defined in Lemma~\ref{lem:3}, and let
$$
\cM_5(x)=\left\{m\in\cM(x)\backslash\(\cup_{i=1}^4\cM_i(x)\):
k\in\cB\right\}.
$$ 
Using~\eqref{eq:newchoiceforP} and Lemma~\ref{lem:3}, we derive that
\begin{eqnarray}
\label{eq:M5}
\#\cM_5(x) & \ll & \frac{x\log_2 x\log_3 x}{(\log x)^2\log_4 x}
\sum_{k\le\cB(x)}\frac {1}{k}\nonumber\\
& \ll & \frac{x\log_3 x}{\log x\log_2 x\log_4 x}
=o\(\frac{x}{\log x}\).
\end{eqnarray}

For integers $m\in\cM(x)\backslash\(\cup_{i=1}^5\cM_i(x)\)$,
the totient $\varphi(k)$ is divisible by every prime 
$$
p\le c_1\frac{\log_2 k}{\log_3 k}.
$$
Since $k>\exp\(\sqrt{\log x}\,\)$, we have
$$
c_1\frac{\log_2 k}{\log_3 k}\ge
\frac{c_1}{2}(1+o(1))\frac{\log_2 x}{\log_3 x}.
$$
Thus, if $x$ is sufficiently large,
$p|\varphi(k)$ for all $p\le y_4=c_2(\log_2 x)/(\log_3 x)$,
where $c_2=\min\{c_1/3,1\}$. Since $k$ and $\varphi(k)$ 
are coprime, it follows that $p\nmid k$ for all primes $p\le y_4$.

Now put $y_5=\log_2 x\log_3 x$, and let 
$$
\cM_6(x)=\{m\in \cM(x)\backslash \(\cup_{i=1}^5\cM_i(x)\):
k\in \cA(x,y_5)\}.
$$
Using Lemma~\ref{lem:2} and the estimate~\eqref{eq:newchoiceforP},
we obtain that
\begin{eqnarray}
\label{eq:M6}
\#\cM_6(x)& \ll & \frac{x\log_2 x\log_3 x}{(\log x)^2\log_4 x}
\sum_{k\in \cA(x,y_5)}\frac{1}{k}\nonumber \\
& \ll & 
\frac{x\log_2 x\log_3 x}{y_5\log x\log_4 x\log_2y_5}
=o\(\frac{x}{\log x}\).
\end{eqnarray}

If $m\in\cM(x)\backslash\(\cup_{i=1}^6\cM_i(x)\)$, then $k$ satisfies
$$
\sum_{\substack{p|(2k-\varphi(k))\\p>y_5}}\frac{1}{p}\le 1.
$$
Note that, since $p|\varphi(k)$ for every prime $p\le y_4$, 
and $p\nmid k$ for any such prime, it follows that
$p\nmid (k-\varphi(k)/2)$ for all $p\le y_4$.  Therefore, 
\begin{eqnarray*}
\sum_{p|(k-\varphi(k)/2)}\frac{1}{p}& \le & \sum_{y_4<p\le y_5} 
\frac{1}{p}+1=\log\(\frac{\log y_4}{\log y_5}\)+O(1)\\
&=&\log\(\frac{\log_3x+\log_4x+O(1)}{\log_3x+\log_4x}\)+O(1)\ll 1,
\end{eqnarray*}
which immediately implies that 
\begin{eqnarray}
\label{eq:phistuff}
\frac{1}{\varphi(k-\varphi(k)/2)} & = & \frac{1}{k-\varphi(k)/2}\cdot 
\frac{k-\varphi(k)/2)}{\varphi(k-\varphi(k)/2)}\nonumber\\
& \ll & 
\frac{1}{k}\prod_{p|(2k-\varphi(k))}\(1+\frac{1}{p-1}\)\nonumber\\
&\le&\frac{1}{k}\exp\(\sum_{p|(2k-\varphi(k))}\frac{1}{p}\)
=\frac{\exp(O(1))}{k}\ll\frac{1}{k}.
\end{eqnarray}

Let $\cM_7(x)=\cM(x)\backslash\(\cup_{i=1}^6\cM_i(x)\)$.
Note that, for every $m\in\cM_7(x)$, the integer $k$ lies
in the set $\cC(x;y_4)$ defined in Lemma~\ref{lem:4}.
Using estimates~\eqref{eq:newchoiceforP}
and~\eqref{eq:phistuff}, together with Lemma~\ref{lem:4}, we derive
that 
\begin{eqnarray}
\label{eq:M7}
\#\cM_7(x)& \ll & \frac{x\log_3 x}{(\log x)^2\log_4 x}\sum_{k\in 
\cC(x;y_4)}\frac{1}{k}\nonumber\\
& \ll & \frac{x\log_3 x}{\log x\log_4 x \log y_4}
=o\(\frac{x}{\log x}\).
\end{eqnarray} 
The assertion of the theorem now follows from estimates 
\eqref{eq:M1}, \eqref{eq:M2}, \eqref{eq:M3}, \eqref{eq:M4}, \eqref{eq:M5},
   \eqref{eq:M6}, and \eqref{eq:M7}.
\end{proof}

\begin{cor}
The infinite series 
$$
\sum_{m\in \cN_c}\frac{1}{m}
$$
is divergent.
\end{cor}

\section{Nonaliquots}
\label{sec:nonaliquot}

\begin{thm}
\label{thm:main2}
The inequality
$$
\#\cN_a(x)\ge\frac{x}{48}(1+o(1))
$$
holds as $x\to\infty$.
\end{thm}

\begin{proof}
Let $\cK$ be the set of positive integers $k\equiv 0\pmod{12}$.
Clearly,
\begin{equation}
\label{eq:1}
\#\cK(x)=\frac{x}{12}+O(1)
\end{equation}

We first determine an upper bound for the cardinality 
of $\(\cK\backslash\cN_a\)(x)$. Let
$k\in\(\cK\backslash\cN_a\)(x)$; then there exists
a positive integer $n$ such that 
$$
f_a(n)=\sigma(n)-n=k.
$$
Since $k\in\cK$, it follows that
\begin{equation}
\label{eq:twelve}
n\equiv\sigma(n)\pmod{12}.
\end{equation}

Assume first that $n$ is odd. Then $\sigma(n)$ is odd as well,
and therefore $n$ is a perfect square.  If $n=p^2$ holds
for some prime $p$, then
$$
x\ge k=\sigma(p^2)-p^2=p+1;
$$
hence, the number of such integers $k$ is at most $\pi(x-1)=o(x)$. 
On the other hand, if $n$ is not the square of a prime, then $n$ has
at least four prime factors (counted with multiplicity).
Let $p_1$ be the smallest prime dividing $n$; then 
$p_1\le n^{1/4}$, and therefore 
$$
n^{3/4}\le \frac{n}{p_1}\le \sigma(n)-n=k\le x;
$$
hence, $n\le x^{4/3}$. Since $n$ is a perfect square, 
the number of integers $k$ is at most $x^{2/3}=o(x)$
in this case.

The above arguments show that all but $o(x)$ integers
$k\in\(\cK\backslash\cN_a\)(x)$ satisfy an equation of the form
$$
f_a(n)=\sigma(n)-n=k
$$
for some {\it even\/} positive integer $n$.  For such $k$, we have
$$
\frac{n}{2}\le\sigma(n)-n=k\le x;
$$
that is, $n\le 2x$.  It follows from the work of~\cite{Er2}
(see, for example, the discussion on page 196 of~\cite{EGPS})
that $12|\sigma(n)$ for all but at most $o(x)$ positive integers
$n\le 2x$.  Hence, using~\eqref{eq:twelve}, we see that every
integer $k\in\(\cK\backslash\cN_a\)(x)$, with at most $o(x)$ exceptions,
can be represented in the form $k=f_a(n)$ for some
$n\equiv 0\pmod{12}$.  For such $k$, we have
$$
x\ge k=\sigma(n)-n=n\left(\frac{\sigma(n)}{n}-1\right)\ge 
n\left(\frac{\sigma(12)}{12}-1\right)=\frac{4n}{3},
$$
therefore $n\le\tfrac{3}{4}x$. Since 
$n$ is a multiple of $12$, it follows that
$$
\#\(\cK\backslash\cN_a\)(x)\le\frac{x}{16}(1+o(1)).
$$
Combining this estimate with~\eqref{eq:1}, we derive that 
\begin{eqnarray*}
\#\cN_a(x)
&\ge&\#\(\cK\cap\cN_a\)(x)=
\#\cK(x)-\#\(\cK\backslash \cN_a\)(x)\\
&\ge&\(\frac{x}{12}-\frac{x}{16}\)(1+o(1))
=\frac{x}{48}(1+o(1)),
\end{eqnarray*}
which completes the proof.
\end{proof}

\section{Robbins numbers}
\label{sec:Robbins}

\begin{thm}
\label{thm:main3}
The inequality
$$
\#\cN_r(x)\ge\frac{x}{3}(1+o(1))
$$
holds as $x\rightarrow \infty$.
\end{thm}

\begin{proof}
Let
\begin{eqnarray*}
\cM_1&=&\{2^\alpha k:k\equiv 3\pmod 6
\text{~and~}\alpha\equiv 0\pmod 2\},\\
\cM_2&=&\{2^\alpha k:k\equiv 5\pmod 6
\text{~and~}\alpha\equiv 1\pmod 2\},
\end{eqnarray*}
and let $\cM$ be the (disjoint) union $\cM_1\cup\cM_2$.  It is easy
to see that
$$
\#\cM_1(x)=\frac{2x}{9}(1+o(1))\qquad\text{and}\qquad
\#\cM_2(x)=\frac{x}{9}(1+o(1))
$$
as $x\to\infty$; therefore,
$$
\#\cM(x)=\frac{x}{3}(1+o(1)).
$$ 
Hence, it suffices to show that all but $o(x)$ numbers in $\cM(x)$
also lie in $\cN_r(x)$.

Let $m\in\cM(x)$, and suppose that $f_r(p)=m$ for some odd prime $p$.
If $m=2^\alpha k$ and $p-1=2^\beta w$, where $k$ and $w$ are positive
and odd, then
$$
2^{\beta-1}(w-\varphi(w))=\frac{p-1}{2}-\varphi(p-1)=f_r(p)=m=2^\alpha k.
$$
If $w=1$, then $w-\varphi(w)=0$, and thus $m=0$, which is not possible. 
Hence, $w\ge 3$, which implies that $\varphi(w)$ is even, and
$w-\varphi(w)$ is odd.  We conclude that 
$\beta=\alpha+1$ and $w-\varphi(w)=k$.

Let us first treat the case that $q^2|w$ for some odd prime $q$.
In this case, we have
$$
k=w-\varphi(w)\ge \frac{w}{q},
$$
and therefore $w\le qk\le qm\le qx$.
Since $q^2|w$ and $w|(p-1)$, it follows that 
$p\equiv 1 \pmod {q^2}$. Note that $q^2\le w\le qx$; hence, $q\le x$.
Since
$$
p=2^{\alpha+1}w+1\le 2^{\alpha+1}qk+1=2qm+1\le 3qx,
$$
the number of such primes $p$ is at most $\pi(3qx;q^2,1)$. 
Put $y=\exp\(\sqrt{\log x}\,\)$.  If $q<x/y$, we use 
again the result of Montgomery and Vaughan~\cite{MV} to derive that 
$$
\pi(3qx;q^2,1)\le
\frac{6qx}{\varphi(q^2)\log(3x/q)}<\frac{6x}{q(q-1)\log y}
<\frac{4x}{q\sqrt{\log x}}
$$
(in the last step, we used the fact that $q\ge 3$), while for $q\ge x/y$, 
we have the trivial estimate 
$$
\pi(3qx;q^2,1)\le\frac{3qx}{q^2}=\frac{3x}{q}.
$$
Summing over $q$, we see that the total number of possibilities
for the prime~$p$ is at most
$$
\frac{4x}{\sqrt{\log x}}\sum_{q<x/y}\frac{1}{q}+3x\sum_{x/y\le q\le x}
\frac{1}{q}.
$$
Since 
$$
\sum_{q<x/y}\frac{1}{q}\ll\log_2(x/y)\le\log_2x,
$$
and
\begin{eqnarray*}
\sum_{x/y\le q\le x}\frac{1}{q}
&=&\log_2x-\log_2(x/y)+O\(\frac{1}{\log x}\)\\
&=&\log\(1+\frac{\log y}{\log x-\log y}\)+O\(\frac{1}{\log x}\)
\ll\frac{1}{\sqrt{\log x}},
\end{eqnarray*}
the number of possibilities for $p$ (hence also for $m=f_r(p)$)
is at most 
$$
O\(\frac{x\log_2x}{\sqrt{\log x}}\)=o(x).
$$
Thus, for the remainder of the proof, we can assume that
$w$ is squarefree.

We claim that $3|w$. Indeed, suppose that this is not the case.
As $w$ is squarefree and coprime to $3$, it follows that
$\varphi(w)\not\equiv 2\pmod 3$ (if $q|w$ for some prime
$q\equiv 1\pmod 3$, then $3|(q-1)|\varphi(w)$;
otherwise $q\equiv 2\pmod 3$ for all $q|w$; hence,
$\varphi(w)=\prod_{q|w}(q-1)\equiv 1\pmod 3$).  In the case that
$m\in\cM_1$, we have $p=2^{\alpha+1}w+1\equiv 2w+1\pmod 3$, thus
$w\not\equiv 1\pmod 3$ (otherwise, $p=3$ and $m=0$);
then $w\equiv 2\pmod 3$.  However, since
$\varphi(w)\not\equiv 2\pmod 3$, it follows that $3$ cannot divide
$k=w-\varphi(w)$, which contradicts the fact that
$k\equiv 3\pmod 6$.  Similarly, in the case that $m\in\cM_2$,
we have $p=2^{\alpha+1}w+1\equiv w+1\pmod 3$, thus
$w\not\equiv 2\pmod 3$; then $w\equiv 1\pmod 3$. However, since
$\varphi(w)\not\equiv 2\pmod 3$, it follows that
$k=w-\varphi(w)\equiv 0$ or $1\pmod 3$, which contradicts the fact
that $k\equiv 5\pmod 6$.  These contradictions establish our claim
that $3|w$.

{From} the preceding result, we have
$$
k=w-\varphi(w)\ge\frac{w}{3},
$$
which implies that
$p=2^{\alpha+1}w+1=2^{\alpha+1}\cdot 3k+1\le 6m+1\le 7x$.
As $\pi(7x)\ll x/\log x$, the number of integers $m\in\cM(x)$ such that
$m=f_r(p)$ for some prime~$p$ of this form
is at most $o(x)$, and this completes the proof.
\end{proof}

\section{Remarks}

Flammenkamp and Luca~\cite{FL} have shown that for every
prime $p$ satisfying the properties: 
\begin{itemize}
\item[$(i)$] $p$ is not Mersenne;
\item[$(ii)$] $p$ is Riesel; i.e., $2^np-1$ is not prime for any $n\ge 1$;
\item[$(iii)$] $2p$ is a noncototient;
\end{itemize}
the number $2^{\ell}p$ is a noncototient for every integer $\ell\ge 0$.
Moreover, they showed that the number of primes $p\le x$ satisfying 
$(i)$ and $(ii)$ is $\gg x/\log x$.  Our Theorem~\ref{thm:main1},
shows that for almost every prime $p$ satisfying $(i)$ and $(ii)$,
$2^\ell p$ is a noncototient for every integer $\ell\ge 0$.  In
particular, these results imply that $\cN_c(x)\ge c(1+o(1))x/\log x$
for some constant $c>1/2$.

It would be interesting to see whether our proof of
Theorem~\ref{thm:main1} can be adapted to show that $\#\cN_c(x)\gg x$,
or to obtain results for the set of positive integers $m$
which are not in the image of the function $n-\lambda(n)$,
where $\lambda(\cdot)$ is the Carmichael function.

\end{document}